\documentclass[12pt]{amsart}
\usepackage[T1]{fontenc}
\usepackage[utf8]{inputenc}
\usepackage[english]{babel}
\usepackage{amsthm,amsmath,amssymb,amsfonts}
\usepackage{a4wide}

\usepackage[usenames]{color}

\usepackage{cite}

\theoremstyle{plain}
\newtheorem{theorem}{Theorem}
\newtheorem{lemma}{Lemma}

\newtheorem*{theorem*}{Theorem}
\newtheorem*{lemma*}{Lemma}

\theoremstyle{definition}

\newcommand{\cof}{ \operatorname {cof} }

\newcommand{\loc}{ {\rm loc} }

\newcommand{\cA}{ \mathcal{A} }

\title{A variational approximation scheme for elastodynamic problems using a new class of admissible mappings}
\author{Anastasia Molchanova}
\thanks{
	This work was partially supported by a
	Grant of the Russian Foundation of 	
	the Russian Science Foundation
	(Agreement  No.~16-41-02004).
	}

\address{Sobolev Institute of Mathematics, 4 Acad. Koptyug avenue, Novosibirsk 630090, Russia}
\address{Peoples' Friendship University, 6 Miklukho-Maklaya str., Moscow 117198, Russia}
\email{a.molchanova@math.nsc.ru}


\begin{document}

\begin{abstract}
	We consider a variational approximation scheme for the 3D elastodynamics problem. Our approach uses a new class of admissible mappings that are closed with respect to the space of mappings with finite distortion.
\vspace{\baselineskip} 
\end{abstract}

\maketitle

%================================================================

The motion of a deformable solid body can be written:
\begin{equation*}
	\frac{\partial^2 y}{\partial t^2} = \nabla \cdot S(\nabla y),
\end{equation*}
where
$y \colon \Omega \times [t_0,t_1] \to \mathbb{R}^3$
is the displacement and
$S$ 
is the first Piola--Kirchgoff stress tensor.
This equation can also be written as a system of conservation laws  
for
the deformation gradient
$F = D_x y$
and
the velocity
$v = {\partial_t} y$
\begin{equation} \label{eq:conserv_law}
	\begin{aligned}
		& \frac{\partial}{\partial t} F_{i \alpha}  = \frac{\partial}{\partial x_{\alpha}} v_{i},\\
		& \frac{\partial}{\partial t} v_{i}  = \sum\limits_{\alpha = 1}\limits^{3} \frac{\partial}{\partial x_{\alpha}} S_{i \alpha} (F),
	\end{aligned}
	\quad i, \alpha = 1, 2, 3.
\end{equation}
In the case of hyperelastic materials, the tensor
$S$ 
can be expressed as the gradient of a scalar function
$W(F)$.
This function is called
the stored energy function,
$W\colon \mathbb{M}^{3\times 3} \to [0, \infty)$,
where 
$\mathbb{M}^{3\times 3}$
stands for $({3 \times 3})$-matrices,
i.e.\
$S(F) = \Big[\dfrac{\partial W}{\partial F_{ij}} (F)\Big]$. 
In order to prove the existence of solutions to the system of equations in \ref{eq:conserv_law}, it would normally be necessary that W be convex.
However, this is incompatible with the known physics of elastic materials such as the requirement of frame-indifference, i.e.\ the principle that certain properties of the system are invariant under arbitrary coordinate transformations \cite{ColNol1959}.

This suggests the replacement of the condition of convexity with a weaker condition such as {\it polyconvexity}
(for further details see, e.g.\ \cite{Ball1977, Ball2002, Ciar1988}).
More precisely, we may assume that 
$$
  W(F) = G(F, \cof F, \det F)
$$
holds for some convex function
$G(F,Z,w)$,
where
$\cof F$ 
and 
$\det F$
are the cofactor matrix 
(i.e.\ transposed adjunctive matrix
$\cof F = {\rm adj}\, F^T$)
and 
determinant of the matrix
$F$
respectively.

At present, the question of
the existence of a solution to
the elastostatic problem has been thoroughly studied.
A review of basic works and open problems
can be found, for example, in
\cite{Ball2002}.
Furthermore,
the reader is referred to 
\cite{Daf2010} 
for local existence of the classical solution of the elastodynamical system  
with rank-one convex and polyconvex stored energy functions.
The existence of global weak solutions,
excluding some particular cases \cite{DiPer1983}, 
is still an open problem.
% except some particular cases.
Nevertheless, the existence of a global measure-valued solution was proven in \cite{DemStuTza2001}
using a variational approximation scheme.

The variational approximation method or,
in the nomenclature of E.~De~Giorgi \cite{DeGio1993},
the minimizing movements method,
is a method by which
the limit of a minimizing sequence of iterations to the variational problem 
for an appropriate functional is found.
The technique developed in \cite{DemStuTza2001} uses the variational approximation scheme to establish a link between elastostatics and elastodynamics.
The method is based on the observation that
the solution of the system~\eqref{eq:conserv_law}
meets the 
{\it additional conservation laws} 
(the idea 
was suggested independently by 
P.\,G.~Le~Floch and T.~Qin \cite{Qin1998})
\begin{equation} \label{eq:ad_conservation_laws}
	\begin{aligned}
		\frac{\partial}{\partial t} \det F & 
			= \sum\limits_{i, \alpha=1}\limits^{3}\frac{\partial}{\partial x_{\alpha}}((\cof F)_{i \alpha} v_i), 
			%\label{id:det} 
		\\
		\frac{\partial}{\partial t}(\cof F)_{k \gamma} & 
			= \sum\limits_{i, j, \alpha, \beta = 1}\limits^{3} \frac{\partial}{\partial x_{\alpha}}
			(\epsilon_{ijk} \epsilon_{\alpha \beta \gamma} F_{j \beta} v_i)
			%\label{id:cof}
	\end{aligned}
\end{equation}
where 
$\epsilon_{ijk}$
is the permutation symbol.
Then, introducing new variables
$Z = \cof F$,
$w = \det F$
and setting
$\Xi = (F, Z, w)$,
we can
use \eqref{eq:conserv_law}
together with \eqref{eq:ad_conservation_laws}
to derive the enlarged system
\begin{equation} \label{enlar_syst}
	\begin{aligned}
		\partial_t v_i &  = \sum\limits_{\alpha,A} \partial_{\alpha} 
			\left(\frac{\partial G}{\partial \Xi^A}(\Xi) \frac{\partial \Phi^A}{\partial F_{i\alpha}(F)} \right) 
			= \sum\limits_{\alpha} \partial_{\alpha} (g_{i\alpha} (\Xi; F)),   \\
		\partial_t \Xi^A & = \sum\limits_{i,\alpha} \partial_{\alpha} 
		\left(\frac{\partial \Phi^A}{\partial F_{i \alpha}}(F)v_i \right),
	\end{aligned}
\end{equation}
where 
\begin{equation*}
	g_{i \alpha} (F, Z, w; F^0) = \frac{\partial G}{\partial F_{i \alpha}} 
		+ \sum\limits_{j, k, \beta, \gamma}  
		\frac{\partial G}{\partial Z_{k \gamma}} \epsilon_{i j k} \epsilon_{\alpha \beta \gamma} F^0_{j \beta}
		+ \frac{\partial G}{\partial w} (\cof F^0)_{i \alpha}. \\
\end{equation*}
and 
$\Phi (F) = (F, \cof F, \det F)$.
Further, time discretization gives
the variational problem
\begin{equation*}
  \min\limits_{\mathcal C} \int\limits_{\Omega} \frac{1}{2}(v - v^0)^2 + G(F,Z,w) \, dx,
\end{equation*}
with initial data 
$v^0(x)$,
$F^0(x)$,
$Z^0(x)$,
$w^0(x)$,
and the admissible set involves the ``additional'' constraints
\begin{equation*}
  \begin{aligned}
    \mathcal{C}= 
      \Bigl\{ (v, F, Z, w) 
      \in
      & L_2 (\Omega) \times L_p (\Omega) \times L_q (\Omega) \times L_r(\Omega),
      \: p>4, \: q,\: r \geq 2 :
      \\
      &\frac{1}{h} (F_{i \alpha} - F^0_{i \alpha}) = \partial_\alpha v_i ,
        \\
      &\frac{1}{h} (Z_{k \gamma} - Z^0_{k \gamma})  
        =  \sum\limits_{i, j, \alpha, \beta} \epsilon_{i j k}\partial_\alpha
        \big( \epsilon_{\alpha \beta \gamma} F^0_{j \beta} v_i \big),
        \\
      &\frac{1}{h} (w - w^0)  
        = \sum\limits_{i, \alpha} \partial_\alpha \big((\cof F^0)_{i \alpha} v_i \big) 
    \Bigr\}. 
  \end{aligned}
\end{equation*}

In this article, we consider the variational approximation scheme 
in a new class of admissible mappings,
in function classes stemming from quasiconformal analysis,
and derive the Euler--Lagrange equations 
in the cases of 
smaller regularity
($p \geq 3$),
finite distortion 
($|F|^3 \leq M w$)
and 
nonnegative Jacobian
($w \geq 0$)
requirements.
Recall that a mapping
$f\colon \Omega \to \mathbb{R}^n$ 
is called the {\it mapping with finite distortion},
$f \in FD(\Omega)$,
if
$f \in W^1_{1, \loc} (\Omega)$,
$J(x, f) \geq 0$
almost everywhere (henceforth abbreviated as a.e.) in 
$\Omega$ 
and
\begin{equation*}
	|Df (x)|^n \le K(x) J(x,f) \quad \text{a.e. in } \Omega,
\end{equation*}
where
$1 \leq K(x) < \infty$
a.e.\ in
$\Omega$
(see for example \cite{IwaSve1993}).
We also note that 
the problem of the approximation preserving the constraint 
$\det F > 0$
is still open, 
except for the very special case of radial elastodynamics \cite{MirTza2012}.
This condition on the deformation gradient is necessary to ensure that the mappings representing motion are orientation-preserving i.e.\ that the deformations are physical.

%\section{Notions}

For the sake of simplicity  
we will work with periodic boundary conditions,
i.e.\ the domain 
$\Omega$
is taken to be a three dimensional torus.
Consider the stored energy function
$W \colon \Omega \times \mathbb{M}^3 \to \mathbb{R}$  
with the following properties:
\smallskip\\
{\bf (H1)} {\bf Polyconvexity:} \label{cond:polyconvexity}
	there exists a convex 
	$C^2$-function 
	$G \colon \mathbb{M}^3\times \mathbb{M}^{3} \times \mathbb{R}_{+} \to \mathbb{R}$ 
	such that for all
	$F\in \mathbb{M}^3$,
	$\det F \geq 0$,
	the equality
	\begin{equation*}
		G(F, \cof F, \det F) = W(F)
	\end{equation*}
	holds.
	\smallskip\\
{\bf (H2)} {\bf Coercivity:} \label{cond:coer}
	there are constants
	$C_1 > 0$,
	$C_2 \in \mathbb{R}$,
	$p \geq 3$,
	$q$,
	$r \geq 2$
%	$s > 2$
	such that
	\begin{equation*}\label{neq:coer}
		G(F, Z, w) \geq C_1 \bigg(|F|^p + |Z|^{q} +  w^r 
%		+ \textcolor{red}{\left(\frac{|F|^n}{w}\right)^{s}} 
		\bigg) + C_2.
	\end{equation*}
	\smallskip\\
{\bf (H3)} \label{cond:finI}
	There is a constant
	$c > 0$
	such that
	\begin{equation*}
		G(F, Z, w) \leq c (|F|^p + |Z|^{q} 
		+ w^r + 1).
	\end{equation*}
	\smallskip\\
{\bf (H4)} \label{cond:limI'}
	There is a constant
	$C > 0$
	such that the inequality
	\begin{equation*}
		|\partial_{F} G|^{p'} + |\partial_{Z} G|^{\frac{p p'}{p-p'}}  
		+ |\partial_{w} G|^{\frac{p p'}{p - 2 p'}} 
		\leq C (|F|^p + |Z|^{q} + w^r + 1)
	\end{equation*}
	holds for
	$p' = \frac{p}{p - 1}$
	if 
	$p > 3$
	and 
	$p' < \frac{3}{2}$
	if 
	$p = 3$.

%\section{} \label{sec:Variational_app_scheme}

Then the iteration scheme is constructed by solving
\begin{eqnarray*} 
	\frac{v_i^J - v_i^{J-1}}{h} &  = & \sum\limits_{\alpha, A} \partial_{\alpha} 
		\left(\frac{\partial G}{\partial \Xi^{A}}(\Xi^{J-1}) 
		\frac{\partial \Phi^A}{\partial F_{i\alpha}(F^{J-1})} \right), \\
	\frac{(\Xi^J - \Xi^{J-1})^A}{h} & = & \sum\limits_{i,\alpha} \partial_{\alpha} 
		\left(\frac{\partial \Phi^A}{\partial F_{i \alpha}}(F^{J-1})v_i^J \right).
\end{eqnarray*}

The 
$J$-th iterates are given by 
\begin{equation*}
	(v^J, \Xi^J) = (v^J, F^J, Z^J,w^J) = (S_h)^J (v^0, F^0, Z^0, w^0),
\end{equation*}
where a solution operator
$S_h$
is defined by 
\begin{subequations} \label{enlar_syst_h}
	\begin{eqnarray}  
		\frac{1}{h} (v_i - v^0_i) & =  &\sum\limits_{\alpha} 
			\partial_\alpha g_{i \alpha}(F, Z, w; F^0), \label{enlar_syst_h:1} \\
		\frac{1}{h} (F_{i \alpha} - F^0_{i \alpha}) & = & \partial_\alpha v_i, \\
		\frac{1}{h} (Z_{k \gamma} - Z^0_{k \gamma}) & = & \sum\limits_{i, j, \alpha, \beta}
			\partial_\alpha (\epsilon_{i j k} \epsilon_{\alpha \beta \gamma} F^0_{j \beta} v_i), \\
		\frac{1}{h} (w - w^0) & = & \sum\limits_{i, \alpha} \partial_\alpha 
			((\cof F^0)_{i \alpha} v_i).
	\end{eqnarray}
\end{subequations}

%\section{}

Given 
$M \in L_s (\Omega)$,
$s >2$,
consider the space
$X = L_2 (\Omega) \times L_p (\Omega) 
\times L_{q} (\Omega) \times L_r(\Omega)$
and the set of admissible mappings
\begin{equation}\label{enlar_syst_h_weak}
  \begin{aligned}
		\cA = 
			& \Bigl\{ (v, F, Z, w) 
			\in X, 
			\: I(v, F, Z, w) < \infty, 
			\: |F(x)|^3 \leq M(x) w(x) \text{ a.e.\ in } \Omega,
			\\
			& 
			\: w(x) \geq 0 \text{ a.e.\ in } \Omega,  
			\text{ and for every } \theta \in C_0^\infty (\Omega, \mathbb{R}^3) 
			\\
			&\int\limits_\Omega \theta \frac{1}{h} (F_{i \alpha} - F^0_{i \alpha}) \, dx 
				= -\int\limits_\Omega v_i \partial_\alpha \theta \, dx, 
				\\
			&\int\limits_\Omega \theta \frac{1}{h} (Z_{k \gamma} - Z^0_{k \gamma}) \, dx 
				= -\int\limits_\Omega  \sum\limits_{i, j, \alpha, \beta} \epsilon_{i j k} \epsilon_{\alpha \beta \gamma} F^0_{j \beta} v_i 
				\partial_\alpha \theta \, dx, 
				\\
			&\int\limits_\Omega \theta \frac{1}{h} (w - w^0) \, dx 
				= -\int\limits_\Omega \sum\limits_{i, \alpha} (\cof F^0)_{i \alpha} 
				v_i \partial_\alpha \theta \, dx 
		\Bigr\}, 
	\end{aligned}
\end{equation}

%\subsection{} 

Let the initial data satisfy
$y^0 = y(0) \in  W^1_p(\Omega) \cap FD(\Omega)$,
$v^0 = \partial_t y(0) \in L_2 (\Omega)$,
$F^0 = D y^0 \in L_p (\Omega)$,
$Z^0 = \cof D y^0 \in L_{q} (\Omega)$,
$w^0 = \det D y^0 \in L_r(\Omega)$,
$|F^0 (x)|^3 \leq M(x) w^0 (x)$,
$w^0 (x) \geq 0$ 
a.e.\ in
$\Omega$
and
$$
	\int\limits_\Omega \frac{1}{2} (v^0)^2 
	+ G(F^0, Z^0, w^0) \, dx < \infty.
$$

%\section{}

It is easy to see that the following assertions hold.
\begin{lemma}
	The admissible set
	$\cA$ 
	is nonempty.
\end{lemma}

\begin{lemma}
	The admissible set
	$\cA$
	is invariant with respect to the relations
	\begin{eqnarray*}
		& \sum\limits_{\alpha} \partial_\alpha Z_{i \alpha} = 0,\\
		& \partial_\beta F_{i \alpha} - \partial_\alpha F_{i \beta} = 0.
	\end{eqnarray*}
	In particular, if
	$F^0$ 
	is a differential 
	then so is
	$F$, 
	and, thus, there exists the mapping
	$y \in W^1_p (\Omega)$ 
	such that
	$\partial_\alpha y_i = F_{i \alpha}$.
\end{lemma}

Consider the minimization problem for the functional
\begin{equation} \label{eq:minI}
	I(v, F, Z, w) = \int\limits_{\Omega} \frac{1}{2} |v - v^0|^2 + G(F, Z, w) \, dx.
\end{equation}

\begin{theorem}\label{thm:exist&uniq}
	There exists
	$(v, F, Z, w) \in \cA$
	satisfying
	\begin{equation*}
		I (v, F, Z, w) 
		= \inf \limits_{\cA} I (v', F', Z', w').
	\end{equation*}
	Furthermore,
	if
	$G$ 
	is a strictly convex function
	then
	the minimizer
	$(v, F, Z, w) \in \cA$
	is unique.
\end{theorem}

The proof of Theorem~\ref{thm:exist&uniq} is based on the next theorem.

\begin{theorem}\label{thm:l.c.s.}
	Let
	$\{(v_n, F_n, Z_n, w_n)\}_{n \in \mathbb{N} }\subset \cA$
	and
	$
		S = \sup\limits_{n\in \mathbb{N}} I (v_n, F_n, Z_n, w_n) < \infty
	$. 
	Then there exist
	$(v, F, Z, w)\in X$
	and a subsequence
	$(v_\mu, F_\mu, Z_\mu, w_\mu)$
	such that
	\begin{equation*} \label{condition:convergence}
		\begin{cases}
			v_\mu \rightharpoonup v & \text{ in } L_2 (\Omega),
			\\
			F_\mu \rightharpoonup F & \text{ in }  L_p(\Omega),
			\\
			Z_\mu \rightharpoonup Z & \text{ in }  L_{q}(\Omega),
			\\
			w_\mu \rightharpoonup w & \text{ in }  L_r(\Omega),
		\end{cases}
	\end{equation*}
	Moreover,
	$(v, F, Z, w)\in \cA$
	and
	\begin{equation*}\label{neq:l.c.s.}
		I(v, F, Z, w) 
		\leq \liminf\limits_{n \rightarrow \infty} I (v_n, F_n, Z_n, w_n) = s < \infty.
	\end{equation*}
\end{theorem}

This statement can be proven by applying the techniques and methods of papers~%
\cite{DemStuTza2001, VodMol2015, MolVod2017}.

%\section{Euler--Lagrange Equations}

We will now show that the minimizer of~%
(\ref{eq:minI})
over the admissible set
$\cA$
satisfies the weak form of the system of equations~%
\eqref{enlar_syst_h}.
To derive the Euler--Lagrange equations, 
we assume that
the minimizer
$( v, F, Z, w) $ 
meets
$w(x) \geq \gamma > 0$
a.e.\ in 
$\Omega$.

We fix ``direction''
${\theta = (\theta_1, \theta_2, \theta_3) \in C^\infty_0 (\Omega, \mathbb{R}^3)}$
such that 
$$
	|F  + h D \theta|^n \leq M \big(w + \sum\limits_{i, \alpha} h \cof F^0_{i \alpha} \, \partial_\alpha \theta_i \big)
	\quad \text{and} \quad
	\sum\limits_{i, \alpha} \cof F^0_{i \alpha} \, \partial_\alpha \theta_i \in L_{\infty} (\Omega).
$$
For
$$
	|\varepsilon| \leq \varepsilon^0 = \frac{\gamma}{\|\sum\limits_{i, \alpha} h 
	\cof F^0_{i \alpha} \partial_\alpha \theta_i\|_{L_\infty} + 1}
$$
we set
\begin{equation*}
	\varepsilon(\delta v_i, \delta F_{i \alpha}, \delta Z_{k \gamma}, \delta w) \\
	= \varepsilon \bigg(\theta_i, h \partial_\alpha \theta_i,  
		\sum\limits_{i, j, \alpha, \beta} h \epsilon_{i j k} \epsilon_{\alpha \beta \gamma} 
		F^0_{j \beta} \partial_\alpha  \theta_i, 
		\sum\limits_{i, \alpha} h \partial_\alpha \cof F^0_{i \alpha} \partial_\alpha \theta_i\bigg)
\end{equation*}
and a variation
\begin{equation*}
	\Xi^\varepsilon = (v^\varepsilon, F^\varepsilon, Z^\varepsilon, w^\varepsilon )
	=  ( v, F, Z, w) + 
	\varepsilon(\delta v_i, \delta F_{i \alpha}, \delta Z_{k \gamma}, \delta w).
\end{equation*}

One can readily see that
such variations fulfill the conditions \eqref{enlar_syst_h_weak}
since columns of the matrix 
$\cof F^0$
are divergence-free vector fields.
Additional requirements on
$\theta$
allow us to conclude that 
$\Xi^\varepsilon = (v^\varepsilon, F^\varepsilon, Z^\varepsilon, w^\varepsilon )$ 
belongs to the admissible set
$\cA$.

Furthermore, using the mean value theorem 
and Lebesgue's dominated convergence theorem
we find that
\begin{multline*}
		\lim \limits_{\varepsilon \to 0} 
		\frac{1}{\varepsilon}  (I(v_i + \varepsilon\delta v_i, F_{i \alpha} + \varepsilon\delta F_{i \alpha},
			Z_{k \gamma} + \varepsilon\delta Z_{k \gamma}, w + \varepsilon\delta w) 
		- I(v_i, F_{i \alpha}, Z_{k \gamma}, w) )\\
		 =  \int\limits_{\Omega}  \sum\limits_{i} \theta_i (v_i - v_i^0)
		+  \sum\limits_{i, \alpha} h \partial_\alpha \theta_i \, g_{i \alpha} (F, Z, w; F^0) \, dx 
\end{multline*}
for
$\varepsilon^* \in [0, \varepsilon]$.
The last equality is the weak form of 
\eqref{enlar_syst_h:1}.
To apply the dominated convergence theorem 
we use the integrability properties of 
$g_{i \alpha}$ 
derived from hypothesis
\textbf{(H4)} 
and Young’s inequality~%
\cite{DemStuTza2001}.

In the case
$p>4$, 
following to \cite{DemStuTza2001},
we find the mapping
$y\colon (0, \infty) \times \Omega \to \mathbb{R}^3$,
$y \in W^{1}_{\infty} ([0,\infty];L_2) \cap L_{\infty} ([0,\infty];W^1_p)$
such that
the conditions
$\partial_t y = v$,
$D_x y = F$,
$\cof D_x y = Z$,
$\det D_x y = w$
are fulfilled.
Moreover, 
they satisfy the weak form of
the additional conservation laws \eqref{eq:ad_conservation_laws} 
and
$y$
is the measure-valued solution of the system~(\ref{enlar_syst}).

\bigskip 

\textbf{Acknowledgment.}
The author is grateful to professor Sergey Vodop$'$yanov for suggesting the problem and useful discussions.

\bibliographystyle{plain}
\bibliography{biblio_ed}

\end{document}